\documentclass[10pt,openbib]{article}
\usepackage{amsmath}
\usepackage{amsfonts}
\usepackage{amssymb}
\usepackage{graphicx}
\newcommand{\R}{\mathbb{R}}

\newcommand{\dee}{\mathop{\! \,\mathrm{d} \!}\nolimits}
\newcommand{\comp}{\raisebox{0pt}{$\scriptstyle\circ \, $}}
\newcommand{\setrule}{\, \rule[-4pt]{.5pt}{13pt}\, }

\newcommand{\onehalf}{\mbox{$\frac{\scriptstyle 1}{\scriptstyle 2}$}}

\newcommand{\ttfrac}[2]{\mbox{$\frac{{\scriptstyle #1}}{{\scriptstyle #2}}$}}

\allowdisplaybreaks

\begin{document}
\thispagestyle{empty}
\begin{center}
{\Large \textbf{The geometry of the Kustaanheimo-Stiefel mapping}}
\mbox{}\vspace{.05in} \\ 
\mbox{}\\
Richard Cushman\footnotemark 
\end{center}
\footnotetext{e-mail: r.h.cushman@gmail.com  
} \medskip 

This paper details the geometry of the Kustaanheimo-Stiefel mapping, which regularizes the 
Hamiltonian of the Kepler problem. It leans heavily on the paper \cite{vdmeer} by J.-C. van der Meer. We use the theory of differential spaces, see \cite[chpt VII \S 3]{cushman-bates} or \cite{sniatycki}.

\section{A $\boldsymbol{{\mathbb{T}}^2}$ action on $T\boldsymbol{{\R }^4}$} \medskip 

Let $T{\R}^4= {\R}^8$ be the tangent bundle of ${\R }^4$ with coordinates $(q,p)$ and 
standard symplectic form $\omega = \sum^4_{i=1}\dee q_i \wedge \dee p_i$. 
Let $\langle \, \, , \, \, \rangle $ be the Euclidean inner product on ${\R }^4$. \medskip 

Consider a Hamiltonian action of the $2$-torus ${\mathbb{T}}^2$ on $T{\R }^4$, which is generated by the flow ${\varphi }^{H_2}_t$ of the Hamiltonian vector field $X_{H_2}$ of the harmonic oscillator 
\begin{displaymath}
H_2(q,p) = \onehalf (\langle p,p \rangle + \langle q, q \rangle )
\end{displaymath} 
and the flow ${\varphi }^{\Xi}_s$ of the Hamiltonian vector field $X_{\Xi}$ of the Hamiltonian 
\begin{displaymath}
\Xi (q,p) = q_1p_2-q_2p_1 +q_3p_4 -q_4p_1.
\end{displaymath}
The momentum mapping of this ${\mathbb{T}}^2$ action is 
\begin{equation}
\mathcal{J}: {\R }^8 \rightarrow {\R }_{\ge 0} \times \R: 
(q,p) \mapsto \big( H_2(q,p), \Xi (q,p) \big). 
\label{eq-zero}
\end{equation}

For every $(h,\xi )$ in the image of $\mathcal{J}$ we determine the ${\mathbb{T}}^2$ reduced space ${\mathcal{J}}^{-1}(h,\xi )/{\mathbb{T}}^2$. We do this in stages. 
First we find the reduced space ${\Xi}^{-1}(\xi )/S^1$ of the $S^1$ action 
generated by ${\varphi }^{\Xi}_s$. The algebra of polynomials on ${\R }^8$, 
which are invariant under the action ${\varphi }^{\Xi}_s$, is generated by 
\begin{equation}
\begin{array}{rlrlrl}
{\pi }_1 = & q^2_1 +q^2_2, & {\pi }_2 = & q^2_3+q^2_4, &  {\pi }_3 =& p^2_1 +p^2_2, \\
{\pi }_4 = & p^2_3+p^2_4, &{\pi }_5 = & q_1p_1+q_2p_2, & {\pi }_6 = &q_3p_3 + q_4p_4, 
\\
{\pi }_7 = & q_1p_2 - q_2p_1&  {\pi }_8 =& q_3p_4 -q_4p_3, & {\pi }_9 = &  
q_1q_4 - q_2q_3, \\
{\pi }_{10} = & q_1q_3+q_2q_4, & {\pi }_{11}= &p_1p_4-p_2p_3, & {\pi }_{12}= &p_1p_3+p_2p_4, \\
{\pi }_{13}= & q_1p_4-q_2p_3, &{\pi }_{14} = & q_1p_3+q_2p_4, & 
{\pi }_{15} = & q_4p_1 -q_3p_2, \\
{\pi }_{16} = & q_3p_1 +q_4p_2. & & 
\end{array}
\label{eq-one}
\end{equation}
The image of the orbit map 
\begin{displaymath}
\Pi : {\R }^8 \rightarrow {\R }^{16}: (q,p) \mapsto \big( {\pi }_1(q,p), \ldots , {\pi }_{16}(q,p) \big)
\end{displaymath}
is the orbit space ${\R }^8/S^1$ of the $S^1$ action ${\varphi }^{\Xi}_s$, which is a semialgebraic 
subset of ${\R }^{16}$ that will be explicitly described below. ${\R }^8/S^1$ is a locally compact 
subcartesian differential space with differential structure $C^{\infty}({\R }^8/S^1)$, where 
$f \in C^{\infty}({\R }^8/S^1)$ if and only if ${\Pi }^{\ast }f \in C^{\infty}({\R }^8)^{S^1}$, the space 
of smooth functions on ${\R }^8$ that are invariant under the $S^1$ action ${\varphi }^{\Xi}_s$. \medskip 

We now show that $C^{\infty}({\R }^8/S^1)$ has a Poisson structure. The quadratic polynomials 
${\pi }_1, \ldots , {\pi }_{16}$ generate a Poisson 
subalgebra of the Poisson algebra of quadratic polynomials $\mathcal{Q}$ on 
$T{\R }^4$, using the Poisson bracket ${\{ \, \, , \, \, \} }_{\mathcal{Q}}$ associated to 
the symplectic form $\omega $, that is, for every $f$, $g \in \in \mathcal{Q}$ one has 
$\{ f,g \}_{\mathcal{Q}}(q,p) = \omega (q,p)\big(X_g(q,p), X_f(q,p) \big)$ for every 
$(q,p) \in T{\R }^4$. Since every smooth function on the orbit space ${\R }^8/S^1$ is a 
smooth function of ${\pi }_1, \ldots , {\pi }_{16}$, the structure matrix 
${\mathcal{W}}_{C^{\infty}({\R }^8/S^1)}$ of the Poisson bracket 
\begin{displaymath}
{\{ f, g \}}_{{\R }^8/S^1} = \sum^{16}_{i,j =1} \frac{\partial g}{\partial {\pi }_1} 
\frac{\partial f}{\partial {\pi }_j} {\{ {\pi }_i, {\pi }_j \}}_{\mathcal{Q}}
\end{displaymath}
for $f$, $g \in C^{\infty}({\R }^8/S^1)$ is ${\mathcal{W}}_{\mathcal{Q}}$. \medskip 

Next we find an explicit description of the orbit space ${\R }^8/S^1$ as a semialgebraic variety in 
${\R }^{16}$. For this we need a new set of generators of the algebra of invariant polynomials. Let 
\begin{subequations}
\begin{align}
K_1 & = -({\pi }_{10} + {\pi }_{12}) =-(q_1q_3+q_2q_4+p_1p_3 +p_2p_4) \notag \\
K_2 &= -({\pi }_9 +{\pi }_{11}) = -(q_1q_4-q_2q_3+p_1p_4-p_2p_3) 
\label{eq-twoa} \\
K_3 & = \onehalf ( {\pi }_2+{\pi }_4 -{\pi }_1 - {\pi }_3) = 
\onehalf (q^2_3+q^2_4 +p^2_3+p^2_4 -q^2_1 -q^2_2 -p^2_1-p^2_2) \notag \\
\rule{0pt}{12pt} L_1 & = {\pi }_{15} - {\pi }_{13} = q_4p_1-q_3p_2+q_2p_3-q_1p_4 
\notag \\
L_2 & = {\pi }_{14} - {\pi }_{16} = q_1p_3+q_2p_4-q_3p_1-q_4p_2 
\label{eq-twob} \\
L_3 & = {\pi }_8-{\pi }_7 = q_3p_4-q_4p_3+q_2p_1-q_1p_2 \notag \\
\rule{0pt}{12pt} H_2 & = \onehalf ({\pi }_1+{\pi }_2+{\pi }_3 +{\pi }_4) = 
\onehalf (q^2_1+q^2_2+q^2_3+q^2_4 +p^2_1 +p^2_2+p^2_3+p^2_4)  \notag \\
\mbox{$\Xi $}\hspace{4pt} & = {\pi }_7+ {\pi }_8 = q_1p_2 - q_2p_1 + q_3p_4 - q_4p_3  
\label{eq-twoc} \\
\rule{0pt}{12pt} U_1 & = -({\pi }_5+ {\pi }_6) = -(q_1p_1+q_2p_2+ q_3p_3+ q_4p_4) 
\notag \\
U_2 & = {\pi }_{10}-{\pi }_{12} = q_1q_3+q_2q_4-p_1p_3-p_2p_4 
\label{eq-twod} \\
U_3 & = {\pi }_9-{\pi }_{11} = q_1q_4 - q_2q_3 +p_2p_3 - p_1p_4 \notag \\
U_4 & = \onehalf ({\pi }_1-{\pi }_2+ {\pi }_4 -{\pi }_3) = \onehalf (q^2_1+q^2_2 
-q^2_3 -q^2_4 +p^2_3 +p^2_4-p^2_1-p^2_2) \notag \\  
\rule{0pt}{12pt} V_1 & = \onehalf ({\pi }_1+{\pi }_2-{\pi }_3-{\pi }_4) = 
\onehalf (q^2_1+q^2_2 +q^2_3+q^2_4-p^2_1 -p^2_2-p^2_3-p^2_4) \notag \\
V_2 & = {\pi }_{14} + {\pi }_{16} = q_1p_3 +q_2p_4 +q_3p_1 + q_4p_2 
\label{eq-twoe} \\
V_3 & = {\pi }_{13}+ {\pi }_{15} = q_1p_4 - q_2p_3 +q_4p_1 - q_3p_2 \notag \\
V_4 & = {\pi }_5-{\pi }_6 = q_1p_1 + q_2p_2 -q_3p_3 -q_4p_4. \notag  
\end{align}
\end{subequations}
The map
\begin{displaymath}
{\R}^{16} \rightarrow {\R }^{16}:({\pi }_1, \ldots , {\pi }_{16}) \mapsto 
\big( K, L, H_2, \Xi ; U, V \big)
\end{displaymath}
is linear and invertible with inverse 
\begin{gather}
\begin{array}{lcl}
{\pi }_1 = \onehalf (H_2-K_3+U_4+V_1) & &
{\pi }_2 = \onehalf (H_2 +K_3-U_4+V_1) \\
\rule{0pt}{12pt}{\pi }_3 = \onehalf (H_2 -K_3-U_4 -V_1) & & 
{\pi }_4  = \onehalf (H_2 +K_3 +U_4-V_1) 
\end{array} \notag \\
\begin{array}{lclcl}
{\pi }_5 =  \onehalf (V_4 -U_1) && {\pi }_6 = -\onehalf (U_1+V_4) && 
{\pi }_7 = \onehalf (\Xi -L_3) \\ 
\rule{0pt}{12pt}{\pi }_8 = \onehalf (\Xi +L_3) & & {\pi }_9 = \onehalf (U_3-K_2) && 
{\pi }_{10} = \onehalf (U_2-K_1) \\ 
\rule{0pt}{12pt} {\pi }_{11} = -\onehalf (U_3 - K_2) && {\pi }_{12} = -\onehalf (U_2+K_1) && {\pi }_{13} = \onehalf (V_3 - L_1) \\ 
\rule{0pt}{12pt} {\pi }_{14} = \onehalf (V_2 + L_2) && 
{\pi }_{15} = \onehalf ( V_3 + L_1) && {\pi }_{16} = \onehalf (V_2 - L_2).
\end{array} \notag
\end{gather}
So $(K,L, H_2, \Xi ;U,V)$ is another set of generators of 
the algebra of polynomials on ${\R }^8$, which are invariant under the 
$S^1$ action ${\varphi }^{\Xi}_s$. Instead of giving the structure matrix of 
the Poisson bracket on $C^{\infty}({\R }^8/S^1)$ we list the Poisson vector fields $Y_G$ 
on ${\R }^8/S^1$ induced from the Hamiltonian vector fields $X_G$ on 
$(T{\R }^4,\omega )$, where $G$ is one of the coordinate functions $(K,L, H_2, \Xi ; U,V)$ listed 
in ($3$).
\begin{align}
Y_{K_1} & = -2L_3\frac{\partial }{\partial K_2} + 2L_2 \frac{\partial }{\partial K_3} 
                       -2K_3\frac{\partial }{\partial L_2} + 2K_2\frac{\partial }{\partial L_3} \notag \\
               & \hspace{,25in} -2U_2 \frac{\partial }{\partial U_1} +2U_2 \frac{\partial }{\partial U_2} 
                     -2V_2\frac{\partial }{\partial V_1} + 2V_1 \frac{\partial }{\partial V_2} \notag \\
Y_{K_2} & = 2L_3\frac{\partial }{\partial K_1} - 2L_1 \frac{\partial }{\partial K_3} 
                       +2K_3\frac{\partial }{\partial L_1} - 2K_1\frac{\partial }{\partial L_3} \notag \\
               & \hspace{.35in} -2U_3 \frac{\partial }{\partial U_1} +2U_1 \frac{\partial }{\partial U_3} 
                     -2V_3\frac{\partial }{\partial V_1} + 2V_1 \frac{\partial }{\partial V_3} \notag \\
Y_{K_3} & = -2L_2\frac{\partial }{\partial K_1} + 2L_1 \frac{\partial }{\partial K_2} 
                       -2K_2\frac{\partial }{\partial L_1} + 2K_1\frac{\partial }{\partial L_2} \notag \\
               & \hspace{.35in} -2U_4 \frac{\partial }{\partial U_1} +2U_1 \frac{\partial }{\partial U_4} 
                     -2V_4\frac{\partial }{\partial V_1} + 2V_1 \frac{\partial }{\partial V_4} \notag \\
Y_{L_1} & = -2K_3\frac{\partial }{\partial K_2} + 2K_2 \frac{\partial }{\partial K_3} 
                       -2L_3\frac{\partial }{\partial L_2} + 2L_2\frac{\partial }{\partial L_3} \notag \\
               & \hspace{.35in} -2U_4 \frac{\partial }{\partial U_3} +2U_3 \frac{\partial }{\partial U_4} 
                     -2V_4\frac{\partial }{\partial V_3} + 2V_3 \frac{\partial }{\partial V_4} \notag \\
Y_{L_2} & = 2K_3\frac{\partial }{\partial K_1} - 2K_1 \frac{\partial }{\partial K_3} 
                       +2L_3\frac{\partial }{\partial L_1} - 2L_1\frac{\partial }{\partial L_3} \notag \\
               & \hspace{.35in} +2U_4 \frac{\partial }{\partial U_2} -2U_2 \frac{\partial }{\partial U_4} 
                     +2V_4\frac{\partial }{\partial V_2} - 2V_2 \frac{\partial }{\partial V_4} \notag \\
Y_{L_3} & = -2K_2\frac{\partial }{\partial K_1} + 2K_1 \frac{\partial }{\partial K_2} 
                       -2L_2\frac{\partial }{\partial L_1} + 2L_1\frac{\partial }{\partial L_2} \notag \\
               & \hspace{.35in} -2U_3 \frac{\partial }{\partial U_2} +2U_2 \frac{\partial }{\partial U_3} 
                     -2V_3\frac{\partial }{\partial V_2} + 2V_2 \frac{\partial }{\partial V_3} \notag \\
Y_{H_2} & = 2V_1\frac{\partial }{\partial U_1} + 2V_2\frac{\partial }{\partial U_2} 
                       +2V_3\frac{\partial }{\partial U_3} + 2V_4\frac{\partial }{\partial U_4} \notag \\
               & \hspace{.35in} -2U_1 \frac{\partial }{\partial V_1} -2U_2 \frac{\partial }{\partial V_2} 
                     -2U_3\frac{\partial }{\partial V_3} - 2U_4 \frac{\partial }{\partial V_4} \notag \\
Y_{\Xi } \hspace{3pt}& = 0 \notag \\ 
Y_{U_1} & = 2U_2\frac{\partial }{\partial K_1} + 2U_3 \frac{\partial }{\partial K_2} 
                       +2U_4\frac{\partial }{\partial K_3} - 2V_1\frac{\partial }{\partial H_2} \notag \\
               & \hspace{.35in} +2K_1 \frac{\partial }{\partial U_2} +2K_2 \frac{\partial }{\partial U_3} 
                     +2K_3\frac{\partial }{\partial U_4} - 2H_2 \frac{\partial }{\partial V_1} \notag \\
Y_{U_2} & = -2U_1\frac{\partial }{\partial K_1} - 2U_4 \frac{\partial }{\partial L_2} 
                       +2U_3\frac{\partial }{\partial L_3} - 2V_2\frac{\partial }{\partial H_2} \notag \\
               & \hspace{.35in} -2K_2 \frac{\partial }{\partial U_1} +2L_1 \frac{\partial }{\partial U_3} 
                     -2L_2\frac{\partial }{\partial U_4} - 2H_2 \frac{\partial }{\partial V_2}  \notag \\ 
Y_{U_3} & = -2U_1\frac{\partial }{\partial K_2} + 2U_4 \frac{\partial }{\partial L_1} 
                       +2U_2\frac{\partial }{\partial L_3} - 2V_3\frac{\partial }{\partial H_2} \notag \\
               & \hspace{.35in} -2K_2 \frac{\partial }{\partial U_1} -2L_1 \frac{\partial }{\partial U_2} 
                     -2V_3\frac{\partial }{\partial U_4} - 2H_2 \frac{\partial }{\partial V_4} \notag \\
Y_{U_4} & = -2U_1\frac{\partial }{\partial K_1} - 2U_3 \frac{\partial }{\partial L_1} 
                       +2U_2\frac{\partial }{\partial L_2} - 2V_4\frac{\partial }{\partial H_2} \notag \\
               & \hspace{.35in} -2K_3 \frac{\partial }{\partial U_1} +2L_2 \frac{\partial }{\partial U_2} 
                     +2V_3\frac{\partial }{\partial U_3} - 2H_2 \frac{\partial }{\partial V_4} \notag \\
Y_{V_1} & = 2V_2\frac{\partial }{\partial K_1} - 2V_4 \frac{\partial }{\partial L_2} 
                       +2V_3\frac{\partial }{\partial L_3} + 2U_2\frac{\partial }{\partial H_2} \notag \\
               & \hspace{.35in} +2H_2 \frac{\partial }{\partial U_2} +2K_1 \frac{\partial }{\partial V_2} 
                     +2K_2\frac{\partial }{\partial V_3} + 2U_4 \frac{\partial }{\partial V_4} \notag \\
Y_{V_2} & = -2V_1\frac{\partial }{\partial K_1} - 2V_4 \frac{\partial }{\partial L_2} 
                       +2V_3\frac{\partial }{\partial L_3} + 2U_2\frac{\partial }{\partial H_2} \notag \\
               & \hspace{.35in} -2H_2 \frac{\partial }{\partial U_2} -2K_1 \frac{\partial }{\partial V_1} 
                     +2L_3\frac{\partial }{\partial V_3} -2L_2 \frac{\partial }{\partial V_4} \notag \\
Y_{V_3} & = -2V_1\frac{\partial }{\partial K_2} + 2V_4 \frac{\partial }{\partial L_1} 
                       -2V_2\frac{\partial }{\partial L_3} + 2U_3\frac{\partial }{\partial H_2} \notag \\
               & \hspace{.35in} +2H_2 \frac{\partial }{\partial U_3} -2K_2 \frac{\partial }{\partial V_1} 
                     -2L_3\frac{\partial }{\partial V_3} +2L_1 \frac{\partial }{\partial V_4} \notag \\
Y_{V_4} & = -2V_1\frac{\partial }{\partial K_3} - 2V_3 \frac{\partial }{\partial L_1} 
                       +2V_2\frac{\partial }{\partial L_2} + 2U_4\frac{\partial }{\partial H_2} \notag \\
               & \hspace{.35in} +2H_2 \frac{\partial }{\partial U_4} -2U_4 \frac{\partial }{\partial V_1} 
                     +2L_2\frac{\partial }{\partial V_2} -2L_1 \frac{\partial }{\partial V_3}. \notag 
\end{align}

\noindent Table 1. \parbox[t]{4.5in}{List of Hamiltonian vector fields of coordinate functions 
on ${\R }^8/S^1$ induced from $(T{\R }^4, \omega )$.} \medskip  

The orbit space ${\R }^8/S^1$ of this $S^1$ action is the $7$ dimensional semialgebraic variety in 
${\R }^{16}$ with coordinates $(K,L, H_2, \Xi ;U,V)$ defined by 
\begin{subequations}
\begin{align} 
\langle U,U \rangle & = U^2_1+U^2_2+U^2_3+U^2_4 = H^2_2 -{\Xi}^2 \ge 0, 
\, \, \, H_2 \ge 0 \notag \\
\langle V,V \rangle & = V^2_1+V^2_2+V^2_3+V^2_4 =  H^2_2 -{\Xi}^2 \ge 0, 
\label{eq-threea} \\
\langle U, V \rangle &= U_1V_1 +U_2V_2+U_3V_3+U_4V_4 = 0, \notag \\
\rule{0pt}{12pt}&\hspace{-.35in}U_2V_1 - U_1V_2  = L_1 \Xi - K_1 H_2, \notag \\
&\hspace{-.35in}U_3V_1 - U_1V_3  = L_2 \Xi - K_2 H_2, \label{eq-threeb} \\
&\hspace{-.35in}U_4V_1 - U_1V_4  = L_3\Xi - K_3 H_2, \notag \\
\rule{0pt}{12pt}&\hspace{-.35in}U_4V_3-U_3V_4  = K_1 \Xi - L_1 H_ 2, \notag \\
&\hspace{-.35in}U_2V_4-U_4V_2  = K_2 \Xi - L_2 H_2, \label{eq-threec} \\
&\hspace{-.35in} U_3V_2-U_2V_3  = K_3\Xi - L_3 H_2. \notag 
\end{align}
\end{subequations}
The orbit map of the $S^1$ action ${\varphi }^{\Xi}_s$ is 
\begin{equation}
\begin{array}{l} \widetilde{\Pi }: {\R }^8 \rightarrow {\R}^8/S^1 \subseteq {\R }^{16}: \\
\hspace{.25in} (q,p) \mapsto \big( K(q,p),L(q,p), H_2(q,p), \Xi (q,p); U(q,p) ,V(q,p) \big),
\end{array} 
\label{eq-four}
\end{equation}
whose image is the orbit space ${\R }^8/S^1$ $(4)$. Since the $S^1$ action generated by 
${\varphi }^{\Xi}_s$ is linear, Schwartz' theorem \cite{schwarz} shows that  
every smooth $S^1$ invariant function on ${\R }^8$ 
is a smooth function of the invariant polynomials $K_j$, $L_j$ for $j=1,2,3$; 
$H_2 $, $\Xi $; and $U_i$, $V_i$ for $i=1, \ldots ,4$. Thus ${\Xi}^{-1}(\xi)/S^1$ is 
defined by setting $\Xi = \xi$ in the relations (\ref{eq-threea}), (\ref{eq-threeb}), and (\ref{eq-threec}). This completes the first stage of studying the ${\mathbb{T}}^2$ action. \medskip 

The second stage in determining ${\mathcal{J}}^{-1}(h, \xi )/{\mathbb{T}}^2$ begins by looking at the harmonic oscillator vector field $X_{H_2} = 
\sum^4_{i=1}(p_i\frac{\partial }{\partial q_i} - q_i\frac{\partial }{\partial p_i})$ on 
$T{\R}^4$. $X_{H_2}$ induces the vector field 
\begin{equation}
Y_{H_2}= 2 \sum^4_{i=1} (V_i \frac{\partial }{\partial U_i} - U_i\frac{\partial }{\partial V_i}) 
\label{eq-five}
\end{equation}
on ${\R }^{16}$, since 
\begin{displaymath}
L_{X_{H_2}}H_2  = L_{X_{H_2}}\Xi = L_{X_{H_2}}K_j = L_{X_{H_2}}L_j =0 \, \, 
\mbox{for $j=1,2,3$}
\end{displaymath}
and 
\begin{displaymath}
L_{X_{H_2}}U_i = 2V_i \quad L_{X_{H_2}}V_i = -2U_i \, \, \mbox{for $i=1,2,3,4$.}
\end{displaymath}
$L_{X_{H_2}}$ leaves invariant the ideal of polynomials, whose 
zeroes define the orbit space ${\R }^8/S^1$, see equation ($4$). The induced vector field $Y_{H_2}$ is a derivation of $C^{\infty}({\R }^8/S^1)$, each of whose integral curves is the image of an integral curve of $X_{H_2}$ on ${\R }^8$ 
under the orbit map ${\Pi }$. The integral curves of $Y_{H_2}$ are defined for all time, 
since $X_{H_2}$ is a complete vector field. Consequently, $Y_{H_2}$ is a complete 
vector field on ${\R }^8/S^1 \subseteq {\R }^{16}$. Its flow is 
the restriction to ${\R }^8/S^1$ of the flow of the vector field $Y_{H_2}$ on ${\R }^{16}$, 
given by
\begin{displaymath}
{\varphi}^{Y_{H_2}}_u(K,L,\Xi ,H_2;U,V) =(0,0,0,0; 2 U \cos u  + 2V \sin u ,  
-2U \sin u + 2V \cos u) . 
\end{displaymath}
This defines a smooth $S^1$ action on the subcartesian space 
${\R }^8/S^1 \subseteq {\R }^{16}$. For $j=1,2,3$ the 
polynomials $K_j$, $L_j$, $\Xi $, and $H_2$ on ${\R }^{16}$ generate the algebra of 
polynomials, which are invariant under the $S^1$ action ${\varphi }^{Y_{H_2}}_u$. Thus the restriction of the preceding polynomials to the orbit space ${\R }^8/S^1$ generate the algebra $C^{\infty}({\R }^8/S^1)^{S^1}$ of smooth functions on ${\R }^8/S^1$, which are 
invariant under the flow ${{\varphi }^{Y_{H_2}}_u}|({\R }^8/S^1)$ of the vector field 
$Y_{H_2}|({\R }^8/S^1)$. \medskip 

We now determine the orbit space of 
the $S^1$ action on ${\R }^8/S^1$ generated by ${{\varphi }^{Y_{H_2}}_u}|({\R }^8/S^1)$. Substituting the equations (\ref{eq-threea}), (\ref{eq-threeb}) and (\ref{eq-threec}) 
into the identity 
\begin{equation}
\sum_{1 \le i < j \le 4}(U_iV_j-U_jV_i)^2 + {\langle U,V\rangle}^2 = 
\langle U, U \rangle \langle V,V \rangle 
\label{eq-six}
\end{equation}
gives the identity
\begin{align*}
(H^2_2 -{\Xi}^2)^2 & = (K^2_1+K^2_2+K^2_3+L^2_1+L^2_2+L^2_3)(H^2_2+{\Xi }^2) \\
& \hspace{1.1in} -4(K_1L_1+K_2L_2+K_3L_3)\Xi H_2, 
\end{align*} 
which holds if 
\begin{subequations}
\begin{align}
0 \le  K^2_1+K^2_2+K^2_3+L^2_1+L^2_2+L^2_3 &= H^2_2+{\Xi }^2, \, \, 
0 \le |\Xi | \le H_2  
\label{eq-sevena} \\
K_1L_1+K_2L_2+K_3L_3 & = \Xi H_2 . \label{eq-sevenb}
\end{align}
\end{subequations}
Using equation ($3$), a calculation shows that equations (\ref{eq-sevena}) and (\ref{eq-sevenb}) indeed hold. Adding and subtracting $\onehalf $ times equation (\ref{eq-sevenb}) from $\ttfrac{1}{4}$ times equation (\ref{eq-sevena}) gives 
\begin{subequations}
\begin{align}
{\xi }^2_1 +{\xi }^2_2 +{\xi }^2_3 & = \ttfrac{1}{4}(H_2 +\Xi )^2, \, \, 
0 \le |\Xi| \le H_2   \label{eq-eighta} \\
{\eta }^2_1 +{\eta }^2_2 +{\eta }^2_3 & = \ttfrac{1}{4}(H_2 - \Xi )^2,  
\label{eq-eightb}
\end{align}
\end{subequations}
where ${\xi }_j = \onehalf (K_j+L_j)$ and ${\eta }_j = \onehalf (K_j - L_j)$ for $j=1,2,3$. 
Equations (\ref{eq-eighta}) and (\ref{eq-eightb}) define the orbit space 
$({\R }^8/S^1)/S^1$ of the $S^1$ action ${{\varphi }^{Y_{H_2}}_u}|({\R }^8/S^1)$ on ${\R }^8/S^1$. The orbit space $({\R }^8/S^1)/S^1$ is the orbit space ${\R }^8/{\mathbb{T}}^2$ of the ${\mathbb{T}}^2$ action on 
${\R }^8$ generated by ${\varphi }^{\Xi }_s$ and ${\varphi }^{H_2}_t$. 
The orbit map of the $S^1$ action ${{\varphi }^{Y_{H_2}}_u}|({\R }^8/S^1)$ on 
${\R }^8/S^1$ induced from the $S^1$ action ${\varphi }^{H_2}_t$ on ${\R }^8$ is 
\begin{equation}
\wp : \big({\R }^8/S^1, C^{\infty}({\R }^8/S^1)\big) \rightarrow \big({\R }^8/{\mathbb{T}}^2, 
C^{\infty}({\R }^8/{\mathbb{T}}^2) \big) , 
\label{eq-nine}
\end{equation}
which is the restriction to ${\R }^8/S^1$ of the smooth map 
\begin{displaymath}
{\R }^{16} \rightarrow {\R}^8 : (K,L,H_2, \Xi ; U,V) \mapsto (\xi , \eta ,H_2, \Xi) = 
\big( \onehalf (K+L), \onehalf (K-L),H_2, \Xi \big) , 
\end{displaymath}
and thus is a smooth mapping of locally compact subcartesian differential spaces. \medskip 

\vspace{-.15in}From equation (\ref{eq-threea}) it follows that the range of the momentum 
map $\mathcal{J}$ (\ref{eq-zero}) of the ${\mathbb{T}}^2$ action on ${\R }^8$ 
is the closed wedge $W = \{ (h, \xi ) \in {\R }_{\ge 0} \times \R \setrule \, 
0 \le |\xi | \le h \}$. When $(h,\xi ) \in \mathrm{int}\, W$, the reduced space 
${\mathcal{J}}^{-1}(h,\xi )/{\mathbb{T}}^2$, defined by equations (\ref{eq-eighta}) and(\ref{eq-eightb}) with $H_2 =h$ and $\Xi = \xi $, is diffeomorphic to $S^2_{\frac{1}{2} (h+\xi )} \times S^2_{\frac{1}{2} (h-\xi )}$. When $(h,\xi ) \in \partial W \setminus \{ (0,0) \}$, that is $\pm \xi = h >0$, the reduced space 
${\mathcal{J}}^{-1}(h,\pm h)/{\mathbb{T}}^2$ is diffeomorphic to $S^2_{h}$. When $(h,\xi ) = (0,0)$ the reduced space ${\mathcal{J}}^{-1}(0,0)/{\mathbb{T}}^2$ is the point 
$(0,0) \in T{\R }^4$. \medskip 

Since the momentum map $\mathcal{J}$ (\ref{eq-zero}) is invariant under the 
$S^1$ action ${\varphi }^{\Xi}_s$, it induces a smooth map
\begin{equation}
J : {\R }^8/S^1 \subseteq {\R }^{16} \rightarrow W \subseteq {\R }_{\ge 0} \times \R:
(K,L, H_2, \Xi ; U,V) \mapsto (H_2, \Xi), 
\label{eq-eleven}
\end{equation}
which is surjective. We now determine the geometry of the fibration defined by $J$. Observe that $J = j \comp \wp$, where 
\begin{displaymath}
j : {\R }^8/{\mathbb{T}}^2 \subseteq {\R }^8 \rightarrow W \subseteq {\R}_{\ge 0} 
\times \R: (\xi ,\eta , H_2, \Xi ) \mapsto (H_2, \Xi )  
\end{displaymath} 
is induced by the momentum mapping of the ${\mathbb{T}}^2$ action on 
$T{\R }^4$. Since $\mathrm{int}\, W$ is simply connected, the fibration 
\begin{displaymath}
J _{|J^{-1}(\mathrm{int}\, W)}: J^{-1}(\mathrm{int}\, W) 
\subseteq {\R }^8/S^1 \rightarrow \mathrm{int}\, W \subseteq {\R }_{\ge 0} \times \R 
\end{displaymath}
is trivial. Hence for each $(h,\xi ) \in \mathrm{int}\, W$ the fiber 
$J^{-1}(h, \xi )$ is diffeomorphic to the fiber $J^{-1}(h,0)$, where 
$(h,0) \in \mathrm{int}\, W$. The fiber $J^{-1}(h,0)$ is defined by 
\begin{gather}
\langle U, U \rangle = h^2  = \langle V, V \rangle , \,  \, \, \langle U,V\rangle =0, \, \, 
h>0 \notag  \\
\begin{array}{l}
-h^{-1}(U_2V_1 - U_1V_2)  = K_1  \\
-h^{-1}(U_3V_1 - U_1V_3)  = K_2  \\
-h^{-1}(U_4V_1 - U_1V_4)  = K_3 
\end{array} \quad 
\begin{array}{l} 
-h^{-1}(U_4V_3 - U_3V_4)  = L_1  \\
-h^{-1}(U_2V_4 - U_4V_2)  = L_2 \\
-h^{-1}(U_3V_2 - U_2V_3)  = L_3  . 
\end{array} 
\label{eq-twelve}
\end{gather}
Hence $J^{-1}(h,0)$ is diffeomorphic to 
\begin{displaymath}
M_{h,0} = \{ (U,V) \in 
{\R }^4\times {\R }^4 \setrule \, \langle U, U \rangle =h^2 = \langle V, V \rangle 
\, \, \& \, \, \langle U,V \rangle =0 \}, 
\end{displaymath}
because $J^{-1}(h,0)$ is the graph of the smooth mapping 
\begin{displaymath}
M_{h,0} \subseteq {\R }^8 \rightarrow {\R }^8: 
(U,V) \mapsto \big( K(U,V), L(U,V),h,0 \big) ,
\end{displaymath}
where $(U,V) \mapsto K(U,V)$ and $(U,V) \mapsto L(U,V)$ are smooth functions 
defined by equation (\ref{eq-twelve}). $M_{h,0}$ is diffeomorphic to $T_hS^3_1$, 
the tangent $h$ sphere bundle to the unit $3$ sphere $S^3_1$, via the mapping 
$(U,V) \mapsto (h^{-1}U, V)$. Since $(\partial W)^{\pm} =\{ (h, \xi ) \in {\R }_{\ge 0} 
\times \R \setrule \, \pm \xi = h >0 \}$ is simply connected, the fibration 
\begin{displaymath}
J_{|J^{-1}((\partial W)^{\pm})}: |J^{-1}((\partial W)^{\pm}) 
\subseteq {\R}^8/S^1\rightarrow (\partial W)^{\pm} \subseteq {\R }_{\ge 0} \times \R
\end{displaymath}
is trivial. Hence for each $(h, \pm h) \in (\partial W)^{\pm}$ with $h>0$ the fiber 
$J^{-1}(h,\pm h)$ is defined by 
\begin{align}
\langle U,U \rangle = \langle V,V \rangle & = h^2, \, \, \langle U,V \rangle =0, \, \, h>0 
\notag \\
h^{-1}(U_2V_1 - U_1V_2) & = \mp h^{-1}(U_4V_3 - U_3V_4) = \pm (L_1-K_1) = \mp {\eta }_1 \notag \\
h^{-1}(U_3V_1 - U_1V_3 & = \mp h^{-1}(U_2V_4 - U_4V_2) = \pm (L_2 -K_2) = \mp {\eta }_2 
\label{eq-thirteen} \\
h^{-1}(U_4V_1 - U_1V_4) & = \mp h^{-1}(U_3V_2 - U_2V_3) = \pm (L_3-K_3) = \mp {\eta }_3. \notag  
\end{align}
Hence $J^{-1}(h, \pm h)$ is diffeomorphic to $M_{h,0}$ because 
it is the graph of the smooth mapping
\begin{displaymath}
M_{h,0} \subseteq {\R }^8 \rightarrow {\R }^8 : 
(U,V) \mapsto \big( \pm (L-K)(U,V), h, \pm h \big) 
\end{displaymath}
where $\pm (L-K)$ is the smooth function defined by equation (\ref{eq-thirteen}). 
When $(h,\xi ) = (0,0) \in W$, the fiber $J^{-1}(0,0)$ is the 
point $(0,0,0,0;0,0) \in {\R}^{16}$, since $h=0$ implies $q = p = 0$, which gives 
$K=L=U=V =\Xi =0$.

\section{The Kustaanheimo-Stiefel mapping}

In this section we define the Kustaanheimo-Stiefel mapping and detail its 
relation to regularizing the Kepler vector field. \medskip 

For $(h,0) \in \mathrm{int}\, W$, the mapping
\begin{displaymath}
\begin{array}{l}
{\wp}_{h,0}: T_hS^3_1 = {\mathcal{J}}^{-1}(h,0)/S^1 \subseteq T{\R }^4 \rightarrow 
S^2_h \times S^2_h = {\mathcal{J}}^{-1}(h,0)/{\mathbb{T}}^2 \subseteq {\R }^6: \\
\rule{0pt}{11pt}(U,V) \mapsto \big( \xi (U,V), \eta (U,V) \big)  = 
\big( \onehalf (K(U,V)+L(U,V)), \onehalf ( K(U,V) - L(U,V)) \big) , 
\end{array} 
\end{displaymath}
where ${\wp}_{h,0}(U,V) = \wp (\xi(U,V),\eta (U,V),h,0;U,V)$, is the orbit map of the preregularized Kepler Hamiltonian (\ref{eq-s2five}). We describe its relation to the Kepler vector field in more detail. \medskip 

On $T_0{\R}^3 = ({\R }^3\setminus \{ 0 \}) \times {\R }^3$ with coordinates 
$(x,y)$ and standard symplectic form $\sum^3_{i=1}\dee x_i \wedge \dee y_i$, 
the Kepler Hamiltonian on the negative energy level set $-\onehalf k^2$ with $k >0$ is 
\begin{equation}
K(x,y) = \onehalf |y|^2 - \frac{1}{|x|} = -\onehalf k^2.
\label{eq-s2one}
\end{equation}
Here $|\, \, |$ is the norm associated to the Euclidean inner product on ${\R }^3$. 
Consider   
\begin{equation}
\widehat{K}(x,y) = \frac{|x|}{k}\big( K(x,y) +\onehalf k^2\big) + \frac{1}{k} = 
\frac{1}{2k}|x| \big( |y|^2 + k^2 \big) . 
\label{eq-s2two}
\end{equation}
The Hamiltonian vector field $X_{\widehat{K}}$ has integral curves which satisfy\
\begin{equation}
\begin{array}{l}
\mbox{\Large $\frac{\dee x}{\dee s} = \frac{|x|}{k} \frac{\partial K}{\partial y}$} \\
\rule{0pt}{16pt} \mbox{\Large $\frac{\dee y}{\dee s} =  
-\frac{|x|}{k}\frac{\partial K}{\partial x}$} -(K(x,y) + \onehalf k^2)\mbox{\Large 
$\frac{\partial \frac{|x|}{k}}{\partial x}$}. 
\end{array}
\label{eq-s2three}
\end{equation} 
On the level set ${\widehat{K}}^{-1}(\frac{1}{k})$ they satisfy 
\begin{equation}
\begin{array}{l}
\mbox{\Large $\frac{\dee x}{\dee s} = \frac{|x|}{k} \frac{\partial K}{\partial y}$} \\
\rule{0pt}{18pt} \mbox{\Large $\frac{\dee y}{\dee s} =  
-\frac{|x|}{k}\frac{\partial K}{\partial x}$.} 
\end{array}
\label{eq-s2four}
\end{equation}
With $\frac{\dee s}{\dee t} = \frac{|x|}{k}$ a solution to equation (\ref{eq-s2four}) 
is a time reparametrization of an integral curve of the Kepler vector field on 
the level set $K^{-1}(-\onehalf k^2)$. The preregularized Kepler Hamiltonian 
on the level set ${\mathcal{K}}^{-1}(1)$ is 
\begin{equation}
\mathcal{K}(x,y) = \onehalf |x|(|y|^2 + 1),  
\label{eq-s2five}
\end{equation}
which is obtained from the Hamiltonian $\widehat{K}$ (\ref{eq-s2two}) using the 
symplectic coordinate change $(x,y) \mapsto (\frac{1}{k}x, ky)$.  \medskip 

On $T_{\ast }{\R }^4 = T{\R }^4 \setminus \{ q = 0 \}$ 
consider the mapping 
\begin{equation}
\mathrm{ks}: T_{\ast }{\R }^4 \rightarrow T_0{\R }^3: 
(q,p) \mapsto (x,y), 
\label{eq-s2fivestar}
\end{equation}
where 
\begin{subequations}
\begin{align}
x_1 & = 2(q_1q_3+q_2q_4) = U_2 -K_1 \notag \\
x_2 &  = 2(q_1q_4- q_2q_3) = U_3 -K_2 
\label{eq-s2sixa} \\
x_3 & = q^2_1+q^2_2-q^2_3 -q^2_4 = U_4-K_3 \notag \\
y_1 & = \frac{1}{\langle q,q \rangle }(q_1p_3+q_2p_4+q_3p_1+q_4p_2) = 
(H_2+V_1)^{-1}V_2 . \notag \\
y_2 & = \frac{1}{\langle q,q \rangle }(q_1p_4-q_2p_3-q_3p_2+q_4p_1) =
(H_2+V_1)^{-1}V_3 
\label{eq-s2sixb} \\
\rule{0pt}{12pt} y_3 & = \frac{1}{\langle q,q \rangle }(q_1p_1+q_2p_2-q_3p_3 -q_4p_4) 
= (H_2+V_1)^{-1}V_4 .\notag 
\end{align}
\end{subequations}
Since the map $\mathrm{ks}$ has components with are invariant under 
${\varphi}^{\Xi}_s$, it sends an orbit of the vector field $X_{\Xi}$ on 
the level set ${\Xi}^{-1}(\xi )\cap T_{\ast }{\R }^4$ to a point in $T_0{\R}^3$. \medskip 

A calculation shows that $\mathrm{ks}$ is 
a Poisson map, that is, $(\mathrm{ks})^{\ast }{\{ f, g \}}_{T{\R }^3} = 
{\{ (\mathrm{ks})^{\ast}f, (\mathrm{ks})^{\ast}g \} }_{T{\R }^4}$ for every 
$f$, $g \in C^{\infty}(T{\R }^3)$. In particular, the map $\mathrm{ks}$ pulls 
back the structure matrix 
\begin{displaymath}%
{\mathcal{W}}_{T{\R }^3}(x,y) = \mbox{\footnotesize $\left( \begin{array}{c|c}
\rule[-5pt]{0pt}{7pt} {\{ x_i,x_j \}}_{T{\R}^3} & {\{ x_i,y_j \}}_{T{\R}^3}  \\ \hline 
\rule{0pt}{10pt} {\{ y_i,x_j \}}_{T{\R}^3} & {\{ y_i,y_j \}}_{T{\R}^3} \end{array} \right) $} 
= \mbox{\footnotesize $\left( \begin{array}{c|c} 0 & 2I_3 \\ \hline -2I_3 & 0 \end{array} 
\right) $}
\end{displaymath}
of the Poisson bracket ${\{ \, , \, \} }_{T{\R }^3}$ on $T{\R }^3$ to the structure matrix
\begin{displaymath}%
{\mathcal{W}}_{T{\R }^4}(q,p) = \mbox{\footnotesize $\left( \begin{array}{c|c}
\rule[-5pt]{0pt}{7pt} {\{ q_i,q_j \}}_{T{\R}^4} & {\{ q_i,p_j \}}_{T{\R}^4}  \\ \hline 
\rule{0pt}{10pt} {\{ p_i,q_j \}}_{T{\R}^4} & {\{ p_i,p_j \}}_{T{\R}^4} \end{array} \right) $} 
= \mbox{\footnotesize $\left( \begin{array}{c|c} 0 & I_4 \\ \hline -I_4 & 0 \end{array} 
\right) $}
\end{displaymath}
of the Poisson bracket ${\{ \, , \, \} }_{T{\R }^4}$ on $T{\R }^4$ associated to 
$\omega = \sum^4_{i=1}\dee q_i \wedge \dee p_i$.  \medskip 

The following calculation determines the pull back of the preregularized \linebreak 
Kepler Hamiltonian $\mathcal{K}$ (\ref{eq-s2five}) by the $\mathrm{ks}$ map on the 
level set ${\mathcal{K}}^{-1}(1)$. Using equation (\ref{eq-s2sixa}) we get
\begin{align*}
|x|^2 & = x^2_1+x^2_2+x^2_3 \\
& = 4(q_2q_3+q_2q_4)^2+4(q_1q_4-q_2q_3)^2 + 
(q^2_1+q^2_2-q^2_3-q^2_4)^2 \\
& = {\langle q,q \rangle }^2 = (H_2+V_1)^2,  \\
|y|^2 & = y^2_1+y^2_2+y^2_3 = \frac{1}{(H_2+V_1)^2}(V^2_1+V^2_2 +V^2_3 
+V^2_4 - V^2_1) \\
& =\frac{1}{(H_2+V_1)^2}(H^2_2-{\Xi }^2 - V^2_1). 
\end{align*}
So 
\begin{align}
(\mathrm{ks})^{\ast }\mathcal{K} & = \onehalf (H_2+V_1) \left[ \frac{1}{(H_2+V_1)^2}(H^2_2 - {\Xi }^2 -V^2_1) \right] 
+ \onehalf (H_2+V_1) \notag \\
& = H_2 - \onehalf \frac{1}{H_2+V_1}{\Xi }^2. 
\label{eq-s2seven} 
\end{align}
Because $\mathrm{ks}$ is a Poisson map, it follows that 
the vector fields $X_{H_2}$ on $T_{\ast }{\R}^4$ and $2X_{\mathcal{K}}$ 
on $T_0{\R}^3$ are $\mathrm{ks}$ related on ${\Xi}^{-1}(0)\cap T_{\ast }{\R }^4$. \medskip

Restricting the mapping $\mathrm{ks}$ (\ref{eq-s2fivestar}) to 
${\Xi }^{-1}(0) \cap T_{\ast }{\R }^4$ gives the Kustaanheimo-Stiefel mapping 
\begin{equation}
\mathrm{KS}: {\Xi }^{-1}(0) \cap T_{\ast }{\R }^4 \rightarrow T_0{\R}^3:(q,p) \mapsto (x,y), 
\label{eq-s2eight}
\end{equation}
where $(x,y)$ are given in equation ($20$). $\mathrm{KS}$ is a Poisson map,  
because the map $\mathrm{ks}$ (\ref{eq-s2fivestar}) is. 
From equation (\ref{eq-s2seven}) it follows that $KS$ pulls back the preregularized Hamiltonian $\mathcal{K}$ (\ref{eq-s2five}) on $T_0{\R }^3$ to the harmonic 
oscillator Hamiltonian $H_2$ on ${\Xi}^{-1}(0)$ and sends an integral curve of the harmonic oscillator vector field $X_{H_2}$ on ${\Xi}^{-1}(0)\cap H^{-1}_2(1)$ to an integral curve of the vector field $2X_{\mathcal{K}}$ on ${\mathcal{K}}^{-1}(1)$, which is just a reparametrization of an integral curve of the Kepler vector field $X_K$ on $K^{-1}(-\onehalf )$. Below we show that those integral curves of $X_{H_2}$ on ${\Xi}^{-1}(0)$ that intersect the plane $\{ q =0 \}$ are sent by the $\mathrm{KS}$ map to integral curves of the Kepler vector field which are collision orbits that reach $\{ x=0 \}$ in finite time, see equation (\ref{eq-ninenwstar}), the $\mathrm{KS}$ map \emph{regularizes} the Kepler vector field on the level set $K^{-1}(-\onehalf )$. 
\medskip 

First we find the pull back of the integrals of angular momentum $J$ and 
eccentricity vector $e$ of the preregularized Kepler vector field on the 
level set ${\mathcal{K}}^{-1}(1)$ by the $\mathrm{KS}$ map and 
determine their relation to the integrals $L_j$ and $K_j$ for $j=1,2,3$ of 
the harmonic oscillator vector field on ${\mathcal{J}}^{-1}(1,0) = H^{-1}_2(1) \cap 
{\Xi }^{-1}(0)$. We use equations (\ref{eq-threea}) and (\ref{eq-threeb}) and the definition of the 
$\mathrm{KS}$ map 
(\ref{eq-s2sixa}) and (\ref{eq-s2sixb}). We start with $L_3$, the third component of 
$L$. On ${\mathcal{J}}^{-1}(1,0)$ we have 
\begin{displaymath}
L_3 = U_2V_3-U_3V_2 = U_2(1+V_1)y_2 -U_3(1+V_1)y_2.
\end{displaymath}
So 
\begin{align*}
(1+V_1)^{-1}L_3 & = (K_1+x_1)y_2 -(K_2+x_2)y_3 = 
K_1y_2-K_2y_3 + x_1y_2 -x_2y_1 \\
& \hspace{-.5in}= (U_1V_2-U_2V_1)\frac{1}{1+V_1}V_3 
-(U_1V_3-U_3V_1)\frac{1}{1+V_1}V_2 +x_1y_2 -x_2y_1 \\
& \hspace{-.5in} = -\frac{1}{1+V_1}V_1(U_2V_3-U_3V_2) +x_1y_2 -x_2y_1 \\
& \hspace{-.5in} = -\frac{1}{1+V_1}V_1L_3 +x_1y_2 -x_2y_1.
\end{align*}
Thus on ${\mathcal{J}}^{-1}(1,0)$ 
\begin{displaymath}
L_3 = [(1+V_1)^{-1}+V_1(1+V_1)^{-1}]L_3 = 
{\mathrm{KS}}^{\ast }( x_1y_2-x_2y_1) = {\mathrm{KS}}^{\ast }(J_3).
\end{displaymath}
A similar argument shows that $L_1$ and $L_2$ on ${\mathcal{J}}^{-1}(1,0)$ 
are equal to the pull back by $\mathrm{KS}$ of $J_1$ and $J_2$, respectively, on 
${\mathcal{K}}^{-1}(1)$. Before dealing with $e_3$, the third component of the eccentricity vector, we compute the pull back of the Euclidean inner product $\langle x,y \rangle $ of $x$ and $y$ on 
${\mathcal{J}}^{-1}(1,0)$ as follows.
\clearpage
\begin{align*}
(\mathrm{KS})^{\ast }(\langle x,y\rangle ) &= (\mathrm{KS})^{\ast }(x_1y_1+x_2y_2+x_3y_3)  \\
&  = (1+V_1)^{-1}\big[ U_2V_2+U_3V_3+U_4V_4 -(V_2K_1+V_3K_2+V_4K_3) \big] \\
&  = -(1+V_1)^{-1}[U_1V_1 + V_2K_1 +V_3K_2 +V_4K_3]\\
&= -(1+V_1)^{-1}[U_1V_1  + (U_1V_2-U_2V_1)V_2 +
(U_1V_3-U_3V_1)V_3  \\
& \hspace{.5in} + (U_1V_4-U_4V_1)V_4] \\
& =-(1+V_1)^{-1}[U_1V_1 +U_1(V^2_1 + V^2_2+V^2_3+V^2_4) \\
&\hspace{.5in}-V_1(U_1V_1+ U_2V_2+U_2V_3+U_4V_4)] \\
& =-(1+V_1)^{-1}[U_1V_1+U_1] = -U_1. 
\end{align*}
Now we compute the pull back by $\mathrm{KS}$ of the third component of the eccentricity vector on ${\mathcal{J}}^{-1}(1,0)$. 
\begin{align*}
(\mathrm{KS})^{\ast }e_3 & = (\mathrm{KS})^{\ast }(-\frac{x_3}{|x|} + 
(y \times (x \times y))_3) 
= (\mathrm{KS})^{\ast }(-\frac{x_3}{|x|} + x_3|y|^2 - y_3 \langle x,y \rangle )\\
& = -(1+V_1)^{-1}[(U_4-K_3) -(U_4-K_3)(1-V_1) +V_4U_1]  \\
& = -(1+V_1)^{-1}[U_4-K_3-U_4+(K_3 -(U_1V_4-U_4V_1))-V_1K_3] \\
& = K_3.
\end{align*}%
A similar argument shows that on ${\mathcal{J}}^{-1}(1,0)$ the pull back 
by $\mathrm{KS}$ of 
first and second component $e_1$ and $e_2$ of the eccentricity vector $e$ are equal to the first and second component $K_1$ and $K_2$, respectively. \medskip 

Next we determine the set $\mathcal{C}$ of initial conditions of the integral curves of $X_{H_2}$ on ${\mathcal{J}}^{-1}(1,0) \subseteq T{\R }^4$, which pass through the collision set $C = \{ q = 0 \}$ in 
$T{\R}^4$. We show that $\mathcal{C}= {\mathcal{J}}^{-1}(1,0) \cap L^{-1}(0)$, that is,   
\begin{equation}
\mathcal{C} = \{ (q,p) \in {\mathcal{J}}^{-1}(1,0) \setrule 
\, L^2_1(q,p) + L^2_2(q,p) +  L^2_3(q,p) = 0 \} . 
\label{eq-ninenwstar}
\end{equation}

\noindent \textbf{Proof.} Suppose that ${\gamma }_{(q,p)}: \R \rightarrow 
{\mathcal{J}}^{-1}(1,0) \subseteq T{\R }^4$ is an integral curve of $X_{H_2}$ 
which starts at $(q,p) \in {\mathcal{J}}^{-1}(1,0)$ and passes through the 
collision set $C$, that is, $(q,p) \in \mathcal{C}$. Then there is a positive time $\tau $ such that 
$q \cos \tau + p \sin \tau =0$, since 
\begin{displaymath}
{\varphi }^{H_2}_t(q,p)=\mbox{\footnotesize $\begin{pmatrix} 
\cos t & \sin t \\ -\sin t & \cos t \end{pmatrix}\begin{pmatrix} q \\ p \end{pmatrix}$.}
\end{displaymath} 
Suppose that $\cos \tau \ne 0$. Then $q = - p \tan \tau = \lambda p$. From the 
definition of the functions $L_j$ for $j=1,2,3$ we get 
\begin{align*}
L_1(\lambda p,p) & = (\lambda p_4)p_1-(\lambda p_2)p_3 + 
(\lambda p_2)p_3 - (\lambda p_1)p_4 = 0 \\
L_2(\lambda p,p) & = (\lambda p_1)p_3 + (\lambda p_2)p_4 - 
(\lambda p_3)p_2 - (\lambda p_4)p_2 = 0 \\
L_3(\lambda p,p) & = (\lambda p_3)p_4-(\lambda p_4)p_3 + 
(\lambda p_2)p_1 - (\lambda p_1)p_2 = 0 . 
\end{align*}
Thus $(q,p) \in {\mathcal{J}}^{-1}(1,0) \cap L^{-1}(0)$. If $\cos \tau =0$, then 
$\sin \tau \ne 0$. So $p = -q \cot \tau = \mu q$. Calculating as above shows that 
$L_(q, \mu ) = L_2(q,\mu q) = L_3(q, \mu q) =0$. So $(q,p) \in L^{-1}(0)\cap 
{\mathcal{J}}^{-1}(1,0) = \mathcal{C}$. \medskip 

Conversely, suppose that $(q,p) \in \mathcal{C}$. Then 
$(q,p) \in {\mathcal{J}}^{-1}(1,0) \cap L^{-1}(0)$. Let $(x,y) =\mathrm{KS}(q,p)$. 
Since $(\mathrm{KS})^{\ast }( J|{\mathcal{K}}^{-1}(1) )= 
L|{\mathcal{J}}^{-1}(1,0)$, it follows 
that $x \times y =0$, because $L(q,p) =0$. Let ${\Gamma }_{(x,y)}$ be 
an integral curve of the vector field $2X_K$ of energy $-1$ starting 
at $(x,y)$, whose angular momentum $J$ vanishes. Then ${\Gamma }_{(x,y)} = 
r(t) e(x,y) = r(t) \frac{x}{|x|}$ is a ray such that $r =r(t) >0$ and $r$ satisfies
\begin{equation}
{\dot{r}}^2 - \frac{2}{r} = -1 . 
\label{eq-s2nine}
\end{equation}
From equation (\ref{eq-s2nine}) we see that $0 \le r \le 2$. Starting at 
$r(0) = r_0$ with $0 < r_0 \le 2$ and $\dot{r}(0) = 0$ there is a finite 
positive time ${\tau }_0$ such that $r({\tau }_0) =0$. To see this separate 
variables in equation (\ref{eq-s2nine}) and integrate. Using the 
change of variables $s^2 = \frac{2}{r} -1 $ we get  
\begin{align*}%
{\tau }_0 & = \int^{{\tau}_0}_0 \dee t = \int^0_{r_0} \frac{\dee r}{\sqrt{2 r^{-1}-1 }} \\ 
&= 4 \int^{s_0}_0 \frac{\dee s}{(s^2+1)^2} , \, \, 
\mbox{where $s_0 = \sqrt{2r^{-1}_0 - 1 } \ge 0$} \\
& < 4\int^{\infty}_0 \frac{\dee s}{(s^2+1)^2} = 
4 \int^{\pi/2}_0 {\cos }^2 u \, \dee u, \, \, 
\mbox{using $s = \tan u$} \\
& =  \pi .
\end{align*}
Since $|x| = \langle q ,q \rangle $, the integral curve ${\gamma }_{(q,p)}$ of 
$X_{H_2}$ starting at $(q,p) \in {\mathcal{J}}^{-1}(1,0)\cap L^{-1}(0)$, whose image under 
the $\mathrm{KS}$ mapping is ${\Gamma }_{(x,y)}$, reaches the collision 
set $C$ at the finite positive time ${\tau }_0$. Hence $(q,p) \in \mathcal{C}$. \hfill $\square $ \medskip 

We now determine the structure matrix ${\mathcal{W}}_{{\R }^8/{\mathbb{T}}^2}$ of the Poisson 
bracket ${\{ \, \, , \, \, \} }_{{\R }^8/{\mathbb{T}}^2}$ on the ${\mathbb{T}}^2$ orbit space 
${\R }^8/{\mathbb{T}}^2$. The smooth surjective mapping 
\begin{displaymath}
\wp : {\R }^8/S^1 \rightarrow {\R }^8/{\mathbb{T}}^2: (K,L,H_2, \Xi ; U,V) \mapsto 
(\xi ,\eta , H_2, \Xi ) 
\end{displaymath}
is a Poisson mapping, that is, ${\wp }^{\ast }{\{ f,g \}}_{{\R }^8/{\mathbb{T}}^2} = 
{\{ {\wp }^{\ast }f, {\wp }^{\ast }g \} }_{{\R }^8/S^1}$, for every $f$, $g \in C^{\infty}({\R }^8/{\mathbb{T}}^2)$. 
Thus we need only determine the Poisson brackets $\{ \, \, , , \, \} $ of the 
functions $K_j$ and $L_j$ for $j =1,2,3$ on $(T{\R }^4, \omega )$. A straightforward computation gives 
\begin{displaymath}
\hspace{-.1in}\{ K_i, K_j \}  = 2 \sum^3_{k=1}{\epsilon}_{ijk} L_k , \, \, 
\{ L_i, L_j \}  = 2 \sum^3_{k=1}{\epsilon}_{ijk} L_k,  \, \, 
\{ K_i, L_j \}  = 2 \sum^3_{k=1}{\epsilon}_{ijk} K_k, 
\end{displaymath}
for $i,j =1,2,3$. Using ${\xi }_j = \onehalf (K_j+L_j)$ and ${\eta }_j = \onehalf (K_j - L_j)$ for $j=1,2,3$ the above equations become 
\begin{equation}
\{ {\xi }_i, {\xi }_j \} =  \sum^3_{k=1}{\epsilon}_{ijk}{\xi }_k, \, \, 
\{ {\eta }_i, {\eta }_j \} =-  \sum^3_{k=1}{\epsilon}_{ijk}{\eta }_k, \, \, \mathrm{and} \, \, 
\{ {\xi }_i, {\eta }_j \} =0.
\label{eq-s2eleven}
\end{equation}


\begin{thebibliography}{9}

\bibitem{cushman-bates} R.H. Cushman and L.M. Bates, 
``Global aspects of classical integrable \linebreak 
systems'', second edition, 
Birkh\"{a}user, Basel, 2015. 

\bibitem{vdmeer} J.-C. van der Meer, Reduction and regularization of the 
Kepler problem, \textit{Celestial. Mech. and Dynam. Astronom.} 
\textbf{132} (2021) 31--50.

\bibitem{schwarz} G. Schwarz, Smooth functions invariant under the action of a
compact {L}ie group, \textit{Topology} \textbf{14} (1975) 63--68. 

\bibitem{sniatycki} J. \'{S}niatycki, ``Differential geometry of singular spaces 
and reduction of \linebreak 
symmetry'', Cambridge University Press, Cambridge, UK, 2013. 
\end{thebibliography}
\end{document}